\documentclass[12pt]{amsart}
\numberwithin{equation}{section} 
\usepackage{amsthm}
\newtheorem{thm}{Theorem}[section]
\newtheorem{lem}[thm]{Lemma}
\newtheorem{prop}[thm]{Proposition}
\newtheorem{cor}[thm]{Corollary}
\newtheorem{conj}[thm]{Conjecture}
\theoremstyle{definition}
\newtheorem{defn}[thm]{Definition}

\theoremstyle{remark}
\newtheorem{rem}[thm]{Remark}
\newcommand{\g}{\mathfrak g}
\newcommand{\gz}{{\mathfrak g}_0}
\newcommand{\bQ}{\mathbf{ Q}}

\begin{document}

\title[Bethe Equation at $q=0$]
{The Bethe Equation at $q=0$, The M\"obius
Inversion Formula, and Weight Multiplicities: \\
III. The $X^{(r)}_N$ case}
%

\author{Atsuo Kuniba}
\address{
Institute of Physics, University of Tokyo,
Tokyo 153-8902, Japan}
\email{atsuo@gokutan.c.u-tokyo.ac.jp}
%
\author{Tomoki Nakanishi}
\address{Graduate School of Mathematics,
    Nagoya University, Nagoya 464-8602, Japan}
\email{nakanisi@math.nagoya-u.ac.jp}
%
%
\author{Zengo Tsuboi}
\address{
Graduate School of Mathematical Sciences, University of Tokyo, 
Komaba, Tokyo 153-8914, Japan}
\email{tsuboi@gokutan.c.u-tokyo.ac.jp}
%

\begin{abstract}
It is shown that the numbers of 
off-diagonal solutions to the $U_q(X^{(r)}_N)$ Bethe equation 
at $q=0$ coincide with the coefficients in the 
recently introduced canonical power series solution of the $Q$-system.
Conjecturally the canonical solutions are characters of the
KR (Kirillov-Reshetikhin) modules.
This implies that the numbers of off-diagonal 
solutions agree with the weight multiplicities, which is interpreted as
a formal completeness of the $U_q(X^{(r)}_N)$ Bethe ansatz 
at $q=0$.
\end{abstract}

\maketitle

%
\section{Introduction}

Enumerating the solutions to the Bethe equation 
began with the invention of the Bethe ansatz \cite{Be}, 
where Bethe himself obtained a counting formula
for $sl_2$-invariant Heisenberg chain.
His calculation is based on the string hypothesis 
and has been generalized to higher spins \cite{K1}, $sl_n$ \cite{K2} and 
a general classical simple Lie algebra $X_n$ \cite{KR}.
These works concern the rational Bethe 
equation \cite{OW}, or in other words, 
$U_q(X^{(1)}_n)$ Bethe equation at $q=1$.

On the other hand, a systematic count at $q=0$ started rather 
recently \cite{KN1,KN2}.
The two approaches are contrastive in many respects.
To explain them, recall the general setting where 
integrable Hamiltonians associated with $U_q(X^{(1)}_n)$ 
act on a finite dimensional module called the quantum space.
At $q=1$, the Hamiltonians are invariant and the Bethe vectors 
are singular with respect to the classical subalgebra $X_n$, while 
for $q \neq 1$, such aspects are no longer valid in general.
Consequently, by completeness at $q=1$ (resp. $q=0$) 
we mean that
the number of solutions to the Bethe equation 
coincides with the multiplicity of  irreducible 
$X_n$ modules (resp. weight multiplicities) in the quantum space.

In this paper we study the Bethe equation associated with 
the quantum affine algebra $U_q(X^{(r)}_N)$ \cite{RW} at $q=0$.
By extending the analyses of the nontwisted case \cite{KN1,KN2},  
an explicit formula $R(\nu,N)$ is derived 
for the number of off-diagonal solutions of the string center equation.
Moreover we relate the result to the {\it $Q$-system} for 
$U_q(X^{(r)}_N)$ introduced in \cite{KR,K3,HKOTT}.
It is a (yet conjectural in general) family of 
character identities for the KR modules (Definition \ref{def:kr}).
Our main finding is that $R(\nu,N)$ is identified with the coefficients in 
the canonical solution of the 
$Q$-system obtained in \cite{KNT}.
Under the Kirillov-Reshetikhin conjecture \cite{KR}
(cf. Conjecture \ref{conj:qsys}), it leads to a 
character formula for tensor products of KR modules, which may be 
viewed as a formal completeness at $q=0$.

The outline of the paper is as follows.
In Section \ref{sec:beq0} we study the $U_q(X^{(r)}_N)$ 
Bethe equation at $q=0$.
For a generic string solution, 
the string centers satisfy the key equation 
(\ref{eq:sce1}), which we call the string center equation (SCE).
There is a one-to-one correspondence between the generic string solutions to 
the Bethe equation and the generic solutions to 
the SCE  (Theorem \ref{thm:C}).
We then enumerate the off-diagonal solutions 
of the SCE, and obtain the formula 
$R(\nu,N)$ in Theorem \ref{thm:Rnds}.
In Section \ref{sec:3} we recall the $Q$-system for $U_q(X^{(r)}_N)$.
It corresponds to a special case (called KR-type) of 
a more general system considered in \cite{KNT}.
There, power series solutions are studied, and 
the notion of the canonical solution is introduced
unifying the ideas in \cite{K1,K2,HKOTY,KN2}.
For the $Q$-system in question, we find that 
the coefficients in the canonical solution are described by $R(\nu,N)$, the
number of off-diagonal solutions of the SCE obtained in Section \ref{sec:3}
(Theorem \ref{th:moohantai}).
A consequence of this fact is stated also in the light of the 
Kirillov-Reshetikhin conjecture \cite{KR,C,KNT}.
We note that the canonical solution of the $Q$-system is 
expressed also as a ratio of two power series \cite{KNT}, which
matches the enumeration at $q=1$ \cite{KR} for the nontwisted cases.

In this paper we omit most of the proofs and calculations, which are 
parallel with those in \cite{KN1,KN2,KNT}.

\section{Bethe equation at $q=0$}\label{sec:beq0}
\subsection{Preliminary}\label{subsec:pre}
Let $\g=X_N$ be a finite-dimensional
complex simple Lie algebra of rank $N$.
We fix a Dynkin diagram automorphism $\sigma$ of $\g$
of order $r=1,2,3$.
The affine Lie algebras of type 
$X^{(r)}_{N}=A^{(1)}_{n} (n \ge 1),B^{(1)}_{n} (n \ge 3),
C^{(1)}_{n} (n \ge 2),D^{(1)}_{n} (n \ge 4),
E^{(1)}_{n} (n=6,7,8), F^{(1)}_{4},G^{(1)}_{2},
A^{(2)}_{2n} (n \ge 1), A^{(2)}_{2n-1} (n \ge 2), 
D^{(2)}_{n+1} (n \ge 2), E^{(2)}_{6}$ and $D^{(3)}_{4}$ 
are realized as the canonical central extension of 
the loop algebras based on the pair $(\g,\sigma)$.
{\unitlength=.95pt
\begin{table}
\caption{Dynkin diagrams for $X^{(r)}_N$.
The enumeration of the nodes with
$I_\sigma \cup \{0\} = \{0,1,\ldots, n\}$ 
is specified under or the right side of the nodes.
In addition, the numbers $d_a$ ($a \in I_\sigma $)
 are attached {\em above} the nodes 
if and only if $d_a \neq 1$.}
\label{tab:Dynkin}
\begin{tabular}[t]{rl}
$A_1^{(1)}$:&
\begin{picture}(26,20)(-5,-5)
\multiput( 0,0)(20,0){2}{\circle{6}}
\multiput(2.85,-1)(0,2){2}{\line(1,0){14.3}}
\put(0,-5){\makebox(0,0)[t]{$0$}}
\put(20,-5){\makebox(0,0)[t]{$1$}}
\put( 6, 0){\makebox(0,0){$<$}}
\put(14, 0){\makebox(0,0){$>$}}
\end{picture}
\\
&
\\
\begin{minipage}[b]{4em}
\begin{flushright}
$A_n^{(1)}$:\\$(n \ge 2)$
\end{flushright}
\end{minipage}&
\begin{picture}(106,40)(-5,-5)
\multiput( 0,0)(20,0){2}{\circle{6}}
\multiput(80,0)(20,0){2}{\circle{6}}
\put(50,20){\circle{6}}
\multiput( 3,0)(20,0){2}{\line(1,0){14}}
\multiput(63,0)(20,0){2}{\line(1,0){14}}
\multiput(39,0)(4,0){6}{\line(1,0){2}}
\put(2.78543,1.1142){\line(5,2){44.429}}
\put(52.78543,18.8858){\line(5,-2){44.429}}
\put(0,-5){\makebox(0,0)[t]{$1$}}
\put(20,-5){\makebox(0,0)[t]{$2$}}
\put(80,-5){\makebox(0,0)[t]{$n\!\! -\!\! 1$}}
\put(100,-7){\makebox(0,0)[t]{$n$}}
\put(55,20){\makebox(0,0)[lb]{$0$}}
\end{picture}
\\
&
\\
\begin{minipage}[b]{4em}
\begin{flushright}
$B_n^{(1)}$:\\$(n \ge 3)$
\end{flushright}
\end{minipage}&
\begin{picture}(126,40)(-5,-5)
\multiput( 20,0)(20,0){2}{\circle{6}}
\multiput(100,0)(20,0){2}{\circle{6}}
\put(10,17.3){\circle{6}}
\put(10,-17.3){\circle{6}}
\multiput(23,0)(20,0){2}{\line(1,0){14}}
\multiput(83,0)(20,0){1}{\line(1,0){14}}
\put(18.5,2.6){\line(-3,5){7.2}}
\put(18.5,-2.6){\line(-3,-5){7.1}}
\multiput(102.85,-1)(0,2){2}{\line(1,0){14.3}} 
\multiput(59,0)(4,0){6}{\line(1,0){2}} 
\put(110,0){\makebox(0,0){$>$}}
\put(10,-22.3){\makebox(0,0)[t]{$1$}}
\put(22,-5){\makebox(0,0)[t]{$2$}}
\put(40,-5){\makebox(0,0)[t]{$3$}}
\put(100,-5){\makebox(0,0)[t]{$n\!\! -\!\! 1$}}
\put(120,-7){\makebox(0,0)[t]{$n$}}
\put(6,8){\makebox(0,0)[l]{$0$}}
\put(10,30.5){\makebox(0,0)[t]{$$}}
\put(9,-4){\makebox(0,0)[t]{$2$}}
\put(22,13){\makebox(0,0)[t]{$2$}}
\put(40,13){\makebox(0,0)[t]{$2$}}
\put(100,13){\makebox(0,0)[t]{$2$}}
\end{picture}
\\
&
\\
\begin{minipage}[b]{4em}
\begin{flushright}
$C_n^{(1)}$:\\$(n \ge 2)$
\end{flushright}
\end{minipage}&
\begin{picture}(126,20)(-5,-5)
\multiput( 0,0)(20,0){3}{\circle{6}}
\multiput(100,0)(20,0){2}{\circle{6}}
\multiput(23,0)(20,0){2}{\line(1,0){14}}
\put(83,0){\line(1,0){14}}
\multiput( 2.85,-1)(0,2){2}{\line(1,0){14.3}} 
\multiput(102.85,-1)(0,2){2}{\line(1,0){14.3}} 
\multiput(59,0)(4,0){6}{\line(1,0){2}} 
\put(10,0){\makebox(0,0){$>$}}
\put(110,0){\makebox(0,0){$<$}}
\put(0,-5){\makebox(0,0)[t]{$0$}}
\put(20,-5){\makebox(0,0)[t]{$1$}}
\put(40,-5){\makebox(0,0)[t]{$2$}}
\put(100,-5){\makebox(0,0)[t]{$n\!\! -\!\! 1$}}
\put(120,-7){\makebox(0,0)[t]{$n$}}
\put(0,13){\makebox(0,0)[t]{$$}}
\put(120,13){\makebox(0,0)[t]{$2$}}
\end{picture}
\\
&
\\
\begin{minipage}[b]{4em}
\begin{flushright}
$D_n^{(1)}$:\\$(n \ge 4)$
\end{flushright}
\end{minipage}&
\begin{picture}(106,40)(-5,-5)
\multiput( 0,0)(20,0){2}{\circle{6}}
\multiput(80,0)(20,0){2}{\circle{6}}
\multiput(20,20)(60,0){2}{\circle{6}}
\multiput( 3,0)(20,0){2}{\line(1,0){14}}
\multiput(63,0)(20,0){2}{\line(1,0){14}}
\multiput(39,0)(4,0){6}{\line(1,0){2}}
\multiput(20,3)(60,0){2}{\line(0,1){14}}
\put(0,-5){\makebox(0,0)[t]{$1$}}
\put(20,-5){\makebox(0,0)[t]{$2$}}
\put(80,-5){\makebox(0,0)[t]{$n\!\! - \!\! 2$}}
\put(103,-5){\makebox(0,0)[t]{$n\!\! -\!\! 1$}}
\put(25,20){\makebox(0,0)[l]{$0$}}
\put(85,20){\makebox(0,0)[l]{$n$}}
\end{picture}
\\
&
\\
$E_6^{(1)}$:&
\begin{picture}(86,60)(-5,-5)
\multiput(0,0)(20,0){5}{\circle{6}}
\multiput(40,20)(0,20){2}{\circle{6}}
\multiput(3,0)(20,0){4}{\line(1,0){14}}
\multiput(40, 3)(0,20){2}{\line(0,1){14}}
\put( 0,-5){\makebox(0,0)[t]{$1$}}
\put(20,-5){\makebox(0,0)[t]{$2$}}
\put(40,-5){\makebox(0,0)[t]{$3$}}
\put(60,-5){\makebox(0,0)[t]{$5$}}
\put(80,-5){\makebox(0,0)[t]{$6$}}
\put(45,20){\makebox(0,0)[l]{$4$}}
\put(45,40){\makebox(0,0)[l]{$0$}}
\end{picture}
\\
&
\\
$E_7^{(1)}$:&
\begin{picture}(126,40)(-5,-5)
\multiput(0,0)(20,0){7}{\circle{6}}
\put(60,20){\circle{6}}
\multiput(3,0)(20,0){6}{\line(1,0){14}}
\put(60, 3){\line(0,1){14}}
\put( 0,-5){\makebox(0,0)[t]{$0$}}
\put(20,-5){\makebox(0,0)[t]{$1$}}
\put(40,-5){\makebox(0,0)[t]{$2$}}
\put(60,-5){\makebox(0,0)[t]{$3$}}
\put(80,-5){\makebox(0,0)[t]{$4$}}
\put(100,-5){\makebox(0,0)[t]{$5$}}
\put(120,-5){\makebox(0,0)[t]{$6$}}
\put(65,20){\makebox(0,0)[l]{$7$}}
\end{picture}
\\
&
\\
\end{tabular}
\begin{tabular}[t]{rl}
$E_8^{(1)}$:&
\begin{picture}(146,40)(-5,-5)
\multiput(0,0)(20,0){8}{\circle{6}}
\put(100,20){\circle{6}}
\multiput(3,0)(20,0){7}{\line(1,0){14}}
\put(100, 3){\line(0,1){14}}
\put( 0,-5){\makebox(0,0)[t]{$0$}}
\put(20,-5){\makebox(0,0)[t]{$1$}}
\put(40,-5){\makebox(0,0)[t]{$2$}}
\put(60,-5){\makebox(0,0)[t]{$3$}}
\put(80,-5){\makebox(0,0)[t]{$4$}}
\put(100,-5){\makebox(0,0)[t]{$5$}}
\put(120,-5){\makebox(0,0)[t]{$6$}}
\put(140,-5){\makebox(0,0)[t]{$7$}}
\put(105,20){\makebox(0,0)[l]{$8$}}
\end{picture}
\\
&
\\
$F_4^{(1)}$:&
\begin{picture}(86,20)(-5,-5)
\multiput( 0,0)(20,0){5}{\circle{6}}
\multiput( 3,0)(20,0){2}{\line(1,0){14}}
\multiput(42.85,-1)(0,2){2}{\line(1,0){14.3}} 
\put(63,0){\line(1,0){14}}
\put(50,0){\makebox(0,0){$>$}}
\put( 0,-5){\makebox(0,0)[t]{$0$}}
\put(20,-5){\makebox(0,0)[t]{$1$}}
\put(40,-5){\makebox(0,0)[t]{$2$}}
\put(60,-5){\makebox(0,0)[t]{$3$}}
\put(80,-5){\makebox(0,0)[t]{$4$}}
\put(0,13){\makebox(0,0)[t]{$$}}
\put(20,13){\makebox(0,0)[t]{$2$}}
\put(40,13){\makebox(0,0)[t]{$2$}}
\end{picture}
\\
&
\\
$G_2^{(1)}$:&
\begin{picture}(46,20)(-5,-5)
\multiput( 0,0)(20,0){3}{\circle{6}}
\multiput( 3,0)(20,0){2}{\line(1,0){14}}
\multiput(22.68,-1.5)(0,3){2}{\line(1,0){14.68}}
\put( 0,-5){\makebox(0,0)[t]{$0$}}
\put(20,-5){\makebox(0,0)[t]{$1$}}
\put(40,-5){\makebox(0,0)[t]{$2$}}
\put(30,0){\makebox(0,0){$>$}}
\put(0,13){\makebox(0,0)[t]{$$}}
\put(20,13){\makebox(0,0)[t]{$3$}}
\end{picture}
\\
&
\\
$A^{(2)}_2$:&
\begin{picture}(26,20)(-5,-5)
\multiput( 0,0)(20,0){2}{\circle{6}}
\multiput(2.958,-0.5)(0,1){2}{\line(1,0){14.084}}
\multiput(2.598,-1.5)(0,3){2}{\line(1,0){14.804}}
\put(0,-5){\makebox(0,0)[t]{$0$}}
\put(20,-5){\makebox(0,0)[t]{$1$}}
\put(10,0){\makebox(0,0){$>$}}
\put(0,13){\makebox(0,0)[t]{$$}}
\end{picture}
\\
&
\\
\begin{minipage}[b]{4em}
\begin{flushright}
$A_{2n}^{(2)}$:\\$(n \ge 2)$
\end{flushright}
\end{minipage}&
\begin{picture}(126,20)(-5,-5)
\multiput( 0,0)(20,0){3}{\circle{6}}
\multiput(100,0)(20,0){2}{\circle{6}}
\multiput(23,0)(20,0){2}{\line(1,0){14}}
\put(83,0){\line(1,0){14}}
\multiput( 2.85,-1)(0,2){2}{\line(1,0){14.3}} 
\multiput(102.85,-1)(0,2){2}{\line(1,0){14.3}} 
\multiput(59,0)(4,0){6}{\line(1,0){2}} 
\put(10,0){\makebox(0,0){$>$}}
\put(110,0){\makebox(0,0){$>$}}
\put(0,-5){\makebox(0,0)[t]{$0$}}
\put(20,-5){\makebox(0,0)[t]{$1$}}
\put(40,-5){\makebox(0,0)[t]{$2$}}
\put(100,-5){\makebox(0,0)[t]{$n\!\! -\!\! 1$}}
\put(120,-7){\makebox(0,0)[t]{$n$}}
\put(0,13){\makebox(0,0)[t]{$$}}
\put(20,13){\makebox(0,0)[t]{$2$}}
\put(40,13){\makebox(0,0)[t]{$2$}}
\put(100,13){\makebox(0,0)[t]{$2$}}
\end{picture}
\\
&
\\
\begin{minipage}[b]{4em}
\begin{flushright}
$A_{2n-1}^{(2)}$:\\$(n \ge 3)$
\end{flushright}
\end{minipage}&
\begin{picture}(126,40)(-5,-5)
\multiput( 0,0)(20,0){3}{\circle{6}}
\multiput(100,0)(20,0){2}{\circle{6}}
\put(20,20){\circle{6}}
\multiput( 3,0)(20,0){3}{\line(1,0){14}}
\multiput(83,0)(20,0){1}{\line(1,0){14}}
\put(20,3){\line(0,1){14}}
\multiput(102.85,-1)(0,2){2}{\line(1,0){14.3}} 
\multiput(59,0)(4,0){6}{\line(1,0){2}} 
\put(110,0){\makebox(0,0){$<$}}
\put(0,-5){\makebox(0,0)[t]{$1$}}
\put(20,-5){\makebox(0,0)[t]{$2$}}
\put(40,-5){\makebox(0,0)[t]{$3$}}
\put(100,-5){\makebox(0,0)[t]{$n\!\! -\!\! 1$}}
\put(120,-7){\makebox(0,0)[t]{$n$}}
\put(25,20){\makebox(0,0)[l]{$0$}}

\put(120,13){\makebox(0,0)[t]{$2$}}
\end{picture}
\\
&
\\
\begin{minipage}[b]{4em}
\begin{flushright}
$D_{n+1}^{(2)}$:\\$(n \ge 2)$
\end{flushright}
\end{minipage}&
\begin{picture}(126,20)(-5,-5)
\multiput( 0,0)(20,0){3}{\circle{6}}
\multiput(100,0)(20,0){2}{\circle{6}}
\multiput(23,0)(20,0){2}{\line(1,0){14}}
\put(83,0){\line(1,0){14}}
\multiput( 2.85,-1)(0,2){2}{\line(1,0){14.3}} 
\multiput(102.85,-1)(0,2){2}{\line(1,0){14.3}} 
\multiput(59,0)(4,0){6}{\line(1,0){2}} 
\put(10,0){\makebox(0,0){$<$}}
\put(110,0){\makebox(0,0){$>$}}
\put(0,-5){\makebox(0,0)[t]{$0$}}
\put(20,-5){\makebox(0,0)[t]{$1$}}
\put(40,-5){\makebox(0,0)[t]{$2$}}
\put(100,-5){\makebox(0,0)[t]{$n\!\! -\!\! 1$}}
\put(120,-7){\makebox(0,0)[t]{$n$}}

\put(20,13){\makebox(0,0)[t]{$2$}}
\put(40,13){\makebox(0,0)[t]{$2$}}
\put(100,13){\makebox(0,0)[t]{$2$}}
\end{picture}
\\
&
\\
$E_6^{(2)}$:&
\begin{picture}(86,20)(-5,-5)
\multiput( 0,0)(20,0){5}{\circle{6}}
\multiput( 3,0)(20,0){2}{\line(1,0){14}}
\multiput(42.85,-1)(0,2){2}{\line(1,0){14.3}} 
\put(63,0){\line(1,0){14}}
\put(50,0){\makebox(0,0){$<$}}
\put( 0,-5){\makebox(0,0)[t]{$0$}}
\put(20,-5){\makebox(0,0)[t]{$1$}}
\put(40,-5){\makebox(0,0)[t]{$2$}}
\put(60,-5){\makebox(0,0)[t]{$3$}}
\put(80,-5){\makebox(0,0)[t]{$4$}}

\put(60,13){\makebox(0,0)[t]{$2$}}
\put(80,13){\makebox(0,0)[t]{$2$}}
\end{picture}
\\
&
\\
$D_4^{(3)}$:&
\begin{picture}(46,20)(-5,-5)
\multiput( 0,0)(20,0){3}{\circle{6}}
\multiput( 3,0)(20,0){2}{\line(1,0){14}}
\multiput(22.68,-1.5)(0,3){2}{\line(1,0){14.68}}
\put( 0,-5){\makebox(0,0)[t]{$0$}}
\put(20,-5){\makebox(0,0)[t]{$1$}}
\put(40,-5){\makebox(0,0)[t]{$2$}}
\put(30,0){\makebox(0,0){$<$}}

\put(40,13){\makebox(0,0)[t]{$3$}}
\end{picture}
\\
&
\\
\end{tabular}
\end{table}}
%
Let $\gz$ be the finite-dimensional $\sigma$-invariant subalgebra of $\g$; namely,

\begin{align*}
\begin{tabular}{c|cccccc}
$\g$ &
$X_n$ & $A_{2n}$ & $A_{2n-1}$ & $D_{n+1}$ 
& $E_6$ & $D_4$ \\
$r$& 1& 2&2 & 2& 2&3\\
\hline
$\gz$ &
$X_n$ & $B_n$ & $C_n$ & $B_n$ & $F_4$ & $G_2$ 
\end{tabular}
\end{align*}
Let $A'=(A'_{ij})$ ($i,j\in I$) and 
$A=(A_{ij})$ ($i,j\in I_\sigma$) be
the Cartan matrices of $\g$ and $\gz$, respectively, where
$I_\sigma$ is the set of $\sigma$-orbits of $I$.
We define the numbers $d'_i$, $d_i$, 
$\epsilon'_i$, $\epsilon_i$ ($i\in I$)
as follows:
$d'_i$ ($i\in I$) are coprime positive
integers such that $(d'_iA'_{ij})$ is symmetric;
$d_i$ ($i\in I_{\sigma}$) are coprime positive
integers such that $(d_iA_{ij})$ is symmetric,
and we set $d_i=d_{\pi(i)}$ ($i\in I$), where 
$\pi: I \rightarrow I_\sigma$ is the canonical projection.
$\epsilon'_i=r$ if $\sigma(i)=i$, and 1
otherwise;
$\epsilon_i=2$ if $A'_{i\sigma(i)}<0$, and 1
otherwise.
Let $\kappa_0=2$ if $X^{(r)}_N
= A^{(2)}_{2n}$, and 1 otherwise.
By the definition one has 
$d_i'=d_i$ and $\epsilon'_i=1$ if\ $r=1$;
$d'_i=1$ if $r>1$; $\epsilon_i=1$ if $X^{(r)}_N
\neq A^{(2)}_{2n}$. 

In this paper we let  $\{1,2,\ldots, N\}$ and 
$\{1,2,\ldots, n\}$ label the sets $I$ and $I_\sigma$, respectively, 
and enumerate  the 
nodes of the Dynkin diagram of $X^{(r)}_N$ by 
$I_\sigma \cup \{0\}$ as specified in 
Table \ref{tab:Dynkin}.
The diagrams (and the enumeration of the nodes for $r>1$) 
coincide with TABLE Aff1-3 in \cite{Kac}, except the 
$A^{(2)}_{2n}$ case.
We fix an injection $\iota: I_\sigma \rightarrow I$ such that 
$\pi \circ \iota = \text{id}_{I_\sigma}$ and 
$A_{a b} < 0 \Leftrightarrow A'_{\iota(a) \iota(b)} < 0$ for any 
$a, b \in I_\sigma$.
To be specific, assume that the labeling of the 
nodes for the Dynkin diagram of 
$\g$ are given by dropping the $0$-th ones {}from $X^{(1)}_N$ 
case in Table \ref{tab:Dynkin}.
Then we simply set $\iota(a) = a$ and regard 
$\iota$ as the embedding of the subset 
$\{1,\ldots, n\} \hookrightarrow \{1,\ldots, N\}$.
The symbols $d'_a, \epsilon'_a$ and $A'_{a b}$ 
for $a,b \in I_\sigma = \{1,\ldots, n\}$ should be interpreted accordingly.
One can check
\begin{align*}
\kappa_0 \epsilon'_a d'_a
=\epsilon_a d_a, \\
\sum_{s=1}^r A'_{a \sigma^s(b)}
=\frac{\epsilon'_a}{\epsilon_a}A_{ab}.
\end{align*}
We use the notation:
\begin{equation}
H = \{(a,m) \mid a \in I_\sigma, m \in {\mathbb Z}_{\ge 1} \}.
\label{eq:hdef1}
\end{equation}

Let $U_{q}(X^{(r)}_N)$ be the quantum affine algebra.
The irreducible finite-dimensional $U_{q}(X^{(r)}_N)$-modules are
parameterized by $N$-tuples of polynomials $(P_i(u))_{i \in I}$
({\em Drinfeld polynomials\/}) with unit constant terms  \cite{CP1, CP2}.
They satisfy the relation 
$P_{\sigma(i)}(u) = P_i(\omega^{\epsilon'_i} u)$, where 
$\omega=\exp(2\pi \sqrt{-1}/r)$. 
Thus it is enough to specify $(P_b(u))_{b \in I_\sigma}$.
Following \cite{KNT} we introduce

\begin{defn}\label{def:kr}
For each $(a,m)\in H$ and $\zeta \in {\mathbb C}^\times$,
let $W^{(a)}_m(\zeta)$ be the finite-dimensional
irreducible $U_q(X^{(r)}_N)$-module whose  Drinfeld
polynomials $P_b(u)$ ($b=1,\dots,n$)
are specified as follows:
$P_b(u)=1$ for $b\neq a$,
and
\begin{align*}
P_a(u)=
\prod_{k=1}^m
(1-\zeta q^{\epsilon_ad_a(m+2-2k)}u).
\end{align*}
We call $W^{(a)}_m(\zeta)$ a
{\em KR (Kirillov-Reshetikhin) module}.
\end{defn}
%
\subsection{The $U_q(X^{(r)}_N)$ Bethe equation} \label{subsec:betheeq}

Let 
\begin{align*}
\mathcal{N}=\{\, N=(N^{(a)}_{m})_{(a,m)\in H}
\mid\, \text{$N^{(a)}_{m}\in \mathbb{Z}_{\geq 0}$},\
\sum_{(a,m)\in H} N^{(a)}_{m} < \infty
\, \}.
\end{align*}
Given $\nu=(\nu^{(a)}_{m})\in {\mathcal N}$, we define 
a tensor product module:
\begin{equation}\label{eq:qspace}
\displaystyle W^\nu =
\bigotimes_{(a,m)\in H}
(W^{(a)}_{m}(\zeta^{(a)}_m))^{\otimes \nu^{(a)}_{m}},
\end{equation}
where $\zeta^{(a)}_m \in {\mathbb C}^\times$.
In the context of solvable lattice models \cite{B}, 
one can regard $W^{\nu}$ as the quantum space on which the commuting 
family of transfer matrices act. 
Reshetikhin and Wiegmann \cite{RW} wrote down 
the $U_q(X^{(r)}_N)$ Bethe equation and conjectured its
relevance to the spectrum of those transfer matrices.
In our formulation, it is the simultaneous equation on 
the complex variables $x_{i}^{(a)}$ 
($i \in \{1,2,\dots,M_{a}\}$, $a \in I_\sigma$) having the form:
\begin{equation}\label{eq:bae-rw}
\prod_{s=1}^{r}\prod_{m=1}^\infty 
\left(
\frac{\omega^{s}(x_{i}^{(a)})^{\frac{1}{\epsilon'_{a}}}
q^{m\kappa_0d'_a\delta_{a,\sigma^{s}(a)}}-1}
{\omega^{s}(x_{i}^{(a)})^{\frac{1}{\epsilon'_{a}}}-
q^{m\kappa_0d'_a\delta_{a,\sigma^{s}(a)}}}
\right)^{\nu^{(a)}_m}
= - \prod_{s=1}^{r} \prod_{b\in I_\sigma}\prod_{j=1}^{M_b}
\frac{\omega^{s}(x_{i}^{(a)})^{\frac{1}{\epsilon'_{a}}}
q^{\kappa_0d'_aA'_{a \sigma^{s}(b)}}
-(x_{j}^{(b)})^{\frac{1}{\epsilon'_{b}}}}
{\omega^{s}(x_{i}^{(a)})^{\frac{1}{\epsilon'_{a}}}-
(x_{j}^{(b)})^{\frac{1}{\epsilon'_{b}}}
q^{\kappa_0d'_aA'_{a \sigma^{s}(b)}}}.
\end{equation} 
For the nontwisted case $r=1$, this reduces to eq.(2.3) in \cite{KN2}.
The both sides are actually rational functions of $(x^{(a)}_i)$.
In the sequel we consider a polynomial version of (\ref{eq:bae-rw}) 
specified as follows:
\begin{equation}\label{eq:bae-fool}
 F^{(a)}_{i+}G^{(a)}_{i-}=
 F^{(a)}_{i-}G^{(a)}_{i+},
\end{equation}
\begin{align*}
F^{(a)}_{i+}&=
 \prod_{k=1}^{\infty}
 (x^{(a)}_{i}q^{k\kappa_0\epsilon'_ad'_a}-1)^{\nu^{(a)}_{k}}, \\ 
F^{(a)}_{i-}&=
 \prod_{k=1}^{\infty}
 (x^{(a)}_{i}-q^{k\kappa_0\epsilon'_ad'_a})^{\nu^{(a)}_{k}}, \\ 
G^{(a)}_{i+}&=
 \prod_{b=1}^{\tilde{n}}
 \prod_{j=1}^{M_{b}}
((x^{(a)}_{i})^{\frac{\epsilon'_{ab}}{\epsilon'_{a}}}
q^{\kappa_0\epsilon'_{ab}d'_aA'_{ab}}-
(x^{(b)}_{j})^{\frac{\epsilon'_{ab}}{\epsilon'_{b}}}), \\ 
G^{(a)}_{i-}&=
 \prod_{b=1}^{\tilde{n}}
 \prod_{j=1}^{M_{b}}
((x^{(a)}_{i})^{\frac{\epsilon'_{ab}}{\epsilon'_{a}}}
-
(x^{(b)}_{j})^{\frac{\epsilon'_{ab}}{\epsilon'_{b}}}
q^{\kappa_0\epsilon'_{ab}d'_aA'_{ab}}),
\end{align*}
where 
$\epsilon'_{a b} = \max(\epsilon'_a, \epsilon'_b)$, 
and $\tilde{n} = n$ except for $\tilde{n} = n+1$ for $A^{(2)}_{2n}$.
When $X^{(r)}_N =  A_{2n}^{(2)}$, 
we have set  $x_{j}^{(n+1)}=-x_{j}^{(n)}$ 
and $M_{n+1}=M_{n}$.

\begin{rem}\label{rem:dpol}
Let ${\mathcal P}^{(a)}_m(u)$ denote the $a$-th Drinfeld polynomial 
of the KR module $W^{(a)}_m(1)$.
Then we have 
\begin{equation*}
\frac{F^{(a)}_{i -}}{F^{(a)}_{i +}} = 
\prod_{(a,m) \in H} 
\left(q^{\epsilon_ad_am}\frac{{\mathcal P}^{(a)}_m(q^{-2\epsilon_ad_a}x^{(a)}_i)}
{{\mathcal P}^{(a)}_m(x^{(a)}_i)}
\right)^{\nu^{(a)}_m}.
\end{equation*}
In view of this, we expect without proof that the solutions of (\ref{eq:bae-rw})
determine the spectrum of transfer matrices acting on 
(\ref{eq:qspace}) with the choice $\zeta^{(a)}_m = 1$.
\end{rem}

We consider  a class of solutions
$(x_i^{(a)})$ of (\ref{eq:bae-fool})
such that  $x_i^{(a)}=x_i^{(a)} (q)$ is meromorphic function of $q$
around $q=0$.
For a meromorphic function $f(q)$ around $q=0$, 
let $\mathrm{ord}(f)$
 be the order of the leading power of the
Laurent expansion of $f(q)$ around $q=0$, i.e.,
\[
f(q) = q^{{\mathrm{ord}}(f)}(f^0 + f^1 q + \cdots ),
 \qquad \ f^0\neq 0,
\]
and let
$\tilde{f}(q) :=  f^0 + f^1q + \cdots$
be the normalized  series.
When $f(q)$ is identically zero, we set $\mathrm{ord}(f)=\infty$.
For each $N=(N^{(a)}_{m})\in {\mathcal N}$, we set 
\begin{align}\label{eq:hdef2}
H'=H'(N):=\{\ (a,m)\in H\mid
N^{(a)}_{m} > 0\ \},
\end{align}
where $H$ is defined in (\ref{eq:hdef1}). 
We have $|H^{\prime}| < \infty$.
\begin{defn}\label{def:string}\upshape
Let  $(M_a)_{a=1}^n$ be the one in
the Bethe equation (\ref{eq:bae-fool}),
and let  $N = (N^{(a)}_{m})\in \mathcal{N}$ satisfy
$ \sum_{m=1}^\infty mN^{(a)}_{m}=M_a$.
A meromorphic solution $(x_i^{(a)})$ of (\ref{eq:bae-fool})
around $q=0$ is called a
{\em string solution of pattern
$N$} if \par
(i) ${{\mathrm{ord}}(F^{(a)}_{i+}G^{(a)}_{i-})}<\infty $ for 
 any $(a,i)$. \par
(ii) $(x_i^{(a)})$ can be arranged as $( x_{m\alpha i}^{(a)})$ with
\begin{equation*}
 (a,m)\in H', \quad \alpha=1,
\dots, N^{(a)}_{m}, \quad i=1,\dots,m
\end{equation*}
such that\par
(a)  $d_{m\alpha i}^{(a)}
 := {\mathrm{ord}}(x_{m\alpha i}^{(a)})=(m+1-2i)\kappa_0\epsilon'_ad'_a$.
\par
(b) $z^{(a)}_{m \alpha}:=x_{m\alpha1}^{(a)0}=x_{m\alpha2}^{(a)0}=
\cdots = x_{m\alpha m}^{(a)0}\ (\neq 0)$,
where $ x_{m\alpha i}^{(a)0}$ is the coefficient of
the leading power of $x_{m\alpha i}^{(a)}$.
\par\noindent
For each $(a,m,\alpha)$,
$(x_{m\alpha i}^{(a)})_{i=1}^m$  is called an {\em $m$-string
of color $a$}, and
 $z^{(a)}_{m \alpha}$  is called
the {\em string center} of the $m$-string
$(x_{m\alpha i}^{(a)})_{i=1}^m$.
Thus, $N^{(a)}_{m}$ is the number of the $m$-strings of color $a$.
\end{defn}
For a string solution 
$x^{(a)}_{m\alpha i}(q)
=q^{d^{(a)}_{m\alpha i}}{\tilde x}^{(a)}_{m\alpha i}(q)$ 
of pattern $N$, the Bethe equation (\ref{eq:bae-fool}) reads
\begin{equation}
 F^{(a)}_{m\alpha i+}G^{(a)}_{m\alpha i-}=
 F^{(a)}_{m\alpha i-}G^{(a)}_{m\alpha i+},
 \label{eq:bae-fg}
\end{equation}
\begin{align}
F^{(a)}_{m\alpha i+}&= \prod_{k=1}^{\infty}
 ({\tilde x}^{(a)}_{m\alpha i}
q^{d^{(a)}_{m\alpha i}+k\kappa_0\epsilon'_ad'_a}-1)^{\nu^{(a)}_{k}}, 
\label{eq:fp}\\ 
F^{(a)}_{m\alpha i-}&=
 \prod_{k=1}^{\infty}
 ({\tilde x}^{(a)}_{m\alpha i}
q^{d^{(a)}_{m\alpha i}}-q^{k\kappa_0\epsilon'_ad'_a})^{\nu^{(a)}_{k}}, 
\label{eq:fm}\\ 
G^{(a)}_{m\alpha i+}&= \prod_{b=1}^{\tilde{n}}\prod_{k=1}^{\infty}
 \prod_{\beta =1}^{N^{(b)}_{k}}\prod_{j=1}^{k}
(({\tilde x}^{(a)}_{m\alpha i})^{\frac{\epsilon'_{ab}}{\epsilon'_{a}}}
q^{\frac{\epsilon'_{ab}}{\epsilon'_{a}}d^{(a)}_{m\alpha i}+
\kappa_0\epsilon'_{ab}d'_aA'_{ab}}-
({\tilde x}^{(b)}_{k\beta j})^{\frac{\epsilon'_{ab}}{\epsilon'_{b}}}
q^{\frac{\epsilon'_{ab}}{\epsilon'_{b}}d^{(b)}_{k\beta j}}), 
\label{eq:gp}\\ 
G^{(a)}_{m\alpha i-}&=
 \prod_{b=1}^{\tilde{n}}\prod_{k=1}^{\infty}
 \prod_{\beta =1}^{N^{(b)}_{k}}\prod_{j=1}^{k}
(({\tilde x}^{(a)}_{m\alpha i})^{\frac{\epsilon'_{ab}}{\epsilon'_{a}}}
q^{\frac{\epsilon'_{ab}}{\epsilon'_{a}}d^{(a)}_{m\alpha i}}-
({\tilde x}^{(b)}_{k\beta j})^{\frac{\epsilon'_{ab}}{\epsilon'_{b}}}
q^{\frac{\epsilon'_{ab}}{\epsilon'_{b}}d^{(b)}_{k\beta j}+
\kappa_0\epsilon'_{ab}d'_aA'_{ab}}),
\label{eq:gm}
\end{align}
where for $X^{(r)}_N = A_{2n}^{(2)}$, we have set 
${\tilde x}^{(n+1)}_{k\beta j}=-{\tilde x}^{(n)}_{k\beta j}$, 
$d^{(n+1)}_{k\beta j}=d^{(n)}_{k\beta j}$ and 
$N^{(n+1)}_{k}=N^{(n)}_{k}$.
According to the procedure similar to \cite{KN2}, we 
can take the $q \to 0$ limit of (\ref{eq:bae-fg}) and 
obtain a key equation:
\begin{equation}
1=(-1)^{m} \prod_{i=1}^{m} 
\frac{F^{(a)0}_{m\alpha i+}G^{(a)0}_{m\alpha i-}}
{F^{(a)0}_{m\alpha i-}G^{(a)0}_{m\alpha i+}}.
\label{eq:bae-lim}
\end{equation}
In order to estimate the order of the Bethe equation (\ref{eq:bae-fg}), we introduce 
\begin{equation}\label{eq:doji}
\begin{split}
\xi^{(a)}_{m\alpha i+}&=
\kappa_0\epsilon'_ad'_a\sum_{k=1}^{\infty}\nu^{(a)}_{k}
  \min(m+1-2i+k,0), \\
\xi^{(a)}_{m\alpha i-}&=
\kappa_0\epsilon'_ad'_a\sum_{k=1}^{\infty}\nu^{(a)}_{k}
  \min(m+1-2i,k), \\ 
\eta^{(a)}_{m\alpha i+}&=\kappa_0
 \sum_{b=1}^{\tilde{n}}\sum_{k=1}^{\infty}
 \sum_{\beta =1}^{N^{(b)}_{k}}\sum_{j=1}^{k}
 \epsilon'_{ab}
 \min(d'_a(m+1-2i+A'_{ab}), d'_b(k+1-2j)), \\ 
\eta^{(a)}_{m\alpha i-}&=
\kappa_0
 \sum_{b=1}^{\tilde{n}}\sum_{k=1}^{\infty}
 \sum_{\beta =1}^{N^{(b)}_{k}}\sum_{j=1}^{k}
\epsilon'_{ab}
  \min(d'_a(m+1-2i), d'_b(k+1-2j+A'_{ba})). 
\end{split}
\end{equation}
\begin{defn}\label{def:generic1}\upshape
A string solution $(x_{m\alpha i}^{(a)})$ to (\ref{eq:bae-fg})
is called {\em generic\/} if 
\begin{align}\label{eq:gerdef}
\begin{split}
{\mathrm{ord}}(F_{m\alpha i\pm}^{(a)}) &= \xi_{m\alpha i\pm}^{(a)},\\
{\mathrm{ord}}(G_{m\alpha i+}^{(a)}) &=
\eta_{m\alpha i+}^{(a)} + \zeta_{m\alpha i}^{(a)},
\qquad
{\mathrm{ord}}(G_{m\alpha i-}^{(a)})=
\eta_{m\alpha i-}^{(a)} + \zeta_{m\alpha i+1}^{(a)},
\end{split}
\end{align}
where  $\zeta_{m\alpha i}^{(a)} :=
{\mathrm{ord}}(\tilde{x}_{m\alpha i}^{(a)}
 - \tilde{x}_{m\alpha i-1}^{(a)})$ for $2 \le i \le m$, and 
$\zeta_{m\alpha1}^{(a)}=\zeta_{m\alpha, m+1}^{(a)}=0$.
\end{defn}
Given a quantum space data $\nu \in \mathcal{N}$ and a string pattern
$N \in \mathcal{N}$, we set
\begin{align}
\gamma^{(a)}_{m} &= \gamma^{(a)}_{m}(\nu)
 = \sum_{k=1}^\infty \min(m,k) \nu_k^{(a)},\label{eq:gammadef}\\
P^{(a)}_{m} &= P^{(a)}_{m}(\nu,N) = \gamma^{(a)}_{m} - \sum_{(b,k)\in H}
\frac{A_{ab}}{\epsilon_ad'_b}\min(d'_am,d'_bk)N_k^{(b)},
\label{eq:pdef}\\
\hat{P}^{(a)}_{m} &= \hat{P}^{(a)}_{m}(\nu,N) = \gamma^{(a)}_{m}
 - \sum_{(b,k)\in H}\frac{A'_{ab}}{d'_b} \min(d'_am,d'_bk)
 N_k^{(b)}.
\label{eq:ppdef}
\end{align}
The number   
$\hat{P}^{(a)}_{m}$ will appear only in the RHS of (\ref{eq:sce1}).
\begin{lem}\label{lem:saa}
We have 
\begin{align}
(\xi^{(a)}_{m\alpha i+}&+\eta^{(a)}_{m\alpha i-})-
(\xi^{(a)}_{m\alpha i-}+\eta^{(a)}_{m\alpha i+}) \nonumber \\ 
&=\begin{cases}
-\kappa_0\epsilon'_ad'_a
(P^{(a)}_{m+1-2i}+N^{(a)}_{m+1-2i})
-\kappa_0\Delta^{(a)}_{m+1-2i} & 1 \le i < \frac{m+1}{2} \nonumber \\
0 & i=\frac{m+1}{2} \\
 \kappa_0\epsilon'_ad'_a
 (P^{(a)}_{2i-m-1}+N^{(a)}_{2i-m-1})
 +\kappa_0\Delta^{(a)}_{2i-m-1} & \frac{m+1}{2} <i \le m ,
\end{cases}
\end{align}
where $\Delta^{(a)}_{j}=0$ except for the following nontwisted 
cases: If there is $a^{\prime}$ such that 
$d_a > d_{a'} = 1$ and $A_{a a'} \neq 0$, then 
\begin{equation*}
 \Delta^{(a)}_{j}=
 \begin{cases}
   -N^{(a^{\prime})}_{2j} & d_{a}=2 \\
   -(N^{(a^{\prime})}_{3j-1}+N^{(a^{\prime})}_{3j}
   +N^{(a^{\prime})}_{3j+1}) & d_{a}=3.
 \end{cases}
\end{equation*}
\end{lem}
For a generic string solution, one can determine the order 
$\zeta^{(a)}_{m\alpha i}$
{}from (\ref{eq:bae-fg}), (\ref{eq:gerdef}) and Lemma \ref{lem:saa}.
Requiring that the resulting 
$\zeta^{(a)}_{m\alpha i}$ should be positive and finite
(cf. Definition \ref{def:string}), one has
\begin{prop}\label{prop:generic}
A necessary condition for the existence of a generic string 
solution of pattern N is 
\begin{equation}
\sum_{k=1}^{\min(i-1,m+1-i)}\left\{ d'_a(P^{(a)}_{m+1-2k}+N^{(a)}_{m+1-2k})
+\Delta^{(a)}_{m+1-2k}\right\}>0,
\label{genaric-condition}
\end{equation}
for $(a,m)\in H^{\prime}, 
1 \le \alpha \le N^{(a)}_{m}, 2 \le i \le m$.
\end{prop}
For a generic string solution, (\ref{eq:bae-lim}) 
becomes an equation for the string centers $(z^{(a)}_{m\alpha})$.
We call it the string center equation (SCE).
\begin{prop}\label{prop:A}
Let $(x_{m\alpha i}^{(a)})$ be a generic string solution of pattern $N$.
Then its string centers $(z^{(a)}_{m\alpha})$ satisfy
the following equations
$((a,m)\in H', \ 1 \le \alpha \le N^{(a)}_{m})$:
\begin{gather}
\prod_{(b,k)\in H'} \prod_{\beta = 1}^{N_k^{(b)}}
(z^{(b)}_{k\beta})^{A_{am\alpha,bk\beta}} =
 (-1)^{\hat{P}^{(a)}_{m}+N^{(a)}_{m}+1},
\label{eq:sce1}\\
A_{am\alpha,bk\beta} = 
\delta_{ab}\delta_{m k}\delta_{\alpha \beta}(P^{(a)}_{m}+N^{(a)}_{m}) +
\frac{A_{ba}}{\epsilon_bd'_a}\min(d'_am,d'_bk) - \delta_{ab}\delta_{m k}.
\label{eq:amat}
\end{gather}
\end{prop}
Note that all the quantities in (\ref{eq:pdef}), (\ref{eq:ppdef}) 
and (\ref{eq:amat}) are integers.
As in \cite{KN2}, Proposition \ref{prop:A} is derived
by explicitly evaluating the ratio
(\ref{eq:bae-lim}) by 
\begin{lem}
For $a\in \{1,2,\dots, n\}$ and $b\in \{1,2,\dots, \tilde{n}\}$, 
we have 
\begin{align*}
&\prod_{i=1}^{m}F^{(a)0}_{m\alpha i \epsilon}=
 \begin{cases}
   (-1)^{\gamma^{(a)}_{m}}
   f_{am\alpha}
   & \epsilon =+ \\ 
   (z^{(a)}_{m\alpha})^{\gamma^{(a)}_{m}}
   f_{am\alpha}
   & \epsilon =-,
 \end{cases}\\
&\prod_{i=1}^{m}\prod_{j=1}^{k}
  ((\tilde{x}^{(a)}_{m\alpha i})^{\frac{\epsilon'_{ab}}{\epsilon'_{a}}}
  q^{\frac{\epsilon'_{ab}}{\epsilon'_{a}}d^{(a)}_{m\alpha i}+
  \frac{1}{2}(1+\epsilon)
  \kappa_0\epsilon'_{ab}d'_aA'_{ab}}-
  (\tilde{x}^{(b)}_{k\beta j})^{\frac{\epsilon'_{ab}}{\epsilon'_{b}}}
  q^{\frac{\epsilon'_{ab}}{\epsilon'_{b}}d^{(b)}_{k\beta j}+
  \frac{1}{2}(1-\epsilon)
  \kappa_0\epsilon'_{ab}d'_aA'_{ab}})^{0} \\ 
&=\begin{cases}
  (-(z^{(b)}_{k\beta})^{\frac{\epsilon'_{ab}}{\epsilon'_{b}}})
  ^{A'_{ab}\min(d'_am,d'_bk)/d'_b
  -\delta_{ab}\delta_{mk}}g^{bk\beta}_{amk} & \epsilon =1 \\ 
  (z^{(a)}_{m\alpha})^{\frac{\epsilon'_{ab}}{\epsilon'_{a}}
  A'_{ab}\min(d'_am,d'_bk)/d'_b-\delta_{ab}\delta_{mk}}
  (-1)^{(m-1)\delta_{ab}\delta_{mk}\delta_{\alpha \beta}}
  g^{bk\beta}_{amk} & \epsilon =-1,
 \end{cases}
\end{align*}
for some $f_{am\alpha}$ and 
$g^{bk\beta}_{am\alpha}$, where we have set $z^{(n+1)}_{m \alpha}:=-z^{(n)}_{m \alpha}$.
\end{lem}
The quantities $f_{am\alpha}$ and 
$g^{bk\beta}_{am\alpha}$ depend on the 
string centers $(z^{(a)}_{m\alpha})$, whose explicit formulae are 
available in \cite{KN2} for nontwisted case.
However we do not need them here.
A string solution is generic if and only if 
$f_{am\alpha}\neq 0$ and 
$g^{bk\beta}_{am\alpha}\neq 0$ for any 
$a \in \{1,\ldots, n\}, b \in \{1, \ldots, \tilde{n}\}, 
m, k \in {\mathbb Z}_{\ge 1}, 
1 \le \alpha \le N^{(a)}_m, 1 \le \beta \le N^{(b)}_k$.
These conditions are equivalent to
\begin{equation}\label{eq:manuke}
\begin{split}
&\mbox{$z^{(a)}_{m\alpha} \neq 1$ {if there is} $k \ge 1$ 
{such that} $\nu^{(a)}_k > 0$ {and} $k \in \langle m \rangle$},\\
&\mbox{$(z^{(a)}_{m\alpha})^{\frac{\epsilon'_{ab}}{\epsilon'_{a}}} \neq 
(z^{(b)}_{k\beta})^{\frac{\epsilon'_{ab}}{\epsilon'_{b}}}$
if $(a,m,\alpha) \neq (b,k,\beta)$ and 
$d'_aA'_{ab} \in \{id'_a-jd'_b \mid i \in \langle m \rangle, 
j \in \langle k \rangle \}$},
\end{split}
\end{equation}
where $\langle m \rangle = \{m-1,m-3,\ldots, -m+1\}$.
Apart {}from the exceptional case $(a,m,\alpha) = (b,k,\beta)$, 
the condition (\ref{eq:manuke}) says that the two terms in each factor in  
(\ref{eq:fp}) -- (\ref{eq:gm}) 
possess different leading terms whenever their orders coincide.
\begin{defn}\label{def:scegene}
A solution to the SCE (\ref{eq:sce1}) is called {\em generic\/}
if it satisfies (\ref{eq:manuke}). 
\end{defn}

Let $A$ be the matrix 
with  the entry $A_{am\alpha, bk\beta}$ in (\ref{eq:amat}). 
The main theorem in this subsection is 
\begin{thm}\label{thm:C}
Suppose that $N \in \mathcal{N}$ satisfies the conditions 
(\ref{genaric-condition})\/
and\/ $\det A\neq 0$.
Then, there is a one-to-one correspondence between
generic string solutions of pattern $N$ 
to the Bethe equation (\ref{eq:bae-fg})
and generic solutions to the SCE (\ref{eq:sce1}) of pattern $N$.
\end{thm}
\begin{rem}
Given the Bethe equation (\ref{eq:bae-rw}), 
the choice of $F^{(a)}_{i\pm}$ and $G^{(a)}_{i \pm}$ in 
(\ref{eq:bae-fool}) is not the unique one.
For example one may restrict the $b$-product in 
$G^{(a)}_{i \pm}$ to those satisfying $A'_{ab} \neq 0$.
Such an ambiguity influences Definition \ref{def:string} (i), 
(\ref{eq:fp}) -- (\ref{eq:gm}),
(\ref{eq:doji}), (\ref{eq:manuke}), hence Definition \ref{def:scegene}.
However, the ratio in (\ref{eq:bae-lim}) is left unchanged, 
and all the statements in Lemma \ref{lem:saa}, 
Propositions \ref{prop:generic}, 
\ref{prop:A} and Theorem \ref{thm:C} remain valid.
\end{rem}
\subsection{Counting of off-diagonal solutions to SCE \label{sec:counting}}
For $k \in \mathbb{C}$ and $j \in \mathbb{Z}$, we define the binomial 
coefficient by
\[
 \binom{k}{j}=\frac{
\varGamma(k+1)}{\varGamma(k-j+1)\varGamma(j+1)}.  
\]
For each $\nu$, $N\in \mathcal{N}$, we define the number
$R(\nu,N)$ by 
\begin{align}
R(\nu,N) &=
\left(\det_{(a,m),(b,k) \in {H'}}F_{am,bk}\right)
\prod_{(a,m) \in H'} \frac{1}{N^{(a)}_{m}}
\binom{P^{(a)}_{m} + N^{(a)}_{m} - 1}{N^{(a)}_{m} - 1},
\label{eq:Rdef}\\
\label{eq:fdef}
F_{am,bk} &=\sum_{\beta =1}^{N^{(b)}_{k}}A_{am\alpha,bk\beta}
= \delta_{ab}\delta_{mk}P^{(a)}_{m}
 + \frac{A_{ba}}{\epsilon_bd'_a}\min(d'_am,d'_bk)N_k^{(b)},
\end{align}
for $N \neq 0$.
Here $H'=H'(N)$ and $P^{(a)}_{m} = P^{(a)}_{m}(\nu,N)$
are given by (\ref{eq:hdef2}) and (\ref{eq:pdef}).
For $N=0$, we set $R(\nu,0) = 1$ irrespective of $\nu$. 
It is easy to see that $R(\nu,N)$ is an integer.
\begin{defn}\label{def:off2}\upshape
A solution $(z_{m\alpha}^{(a)})$ to the SCE
is called {\em off-diagonal\/} ({\em diagonal\/}) if
$z^{(a)}_{m\alpha} = z^{(a)}_{m\beta}$ only for $\alpha=\beta$
(otherwise).
\end{defn}
Our main result in this subsection is
\begin{thm}\label{thm:Rnds}
Suppose  $P^{(a)}_m(\nu,N) \ge 0$ 
for any $(a,m) \in H'$.
Then the number of off-diagonal solutions to
the SCE (\ref{eq:sce1}) of pattern $N$
divided by $\prod_{(a,m)\in H'} N^{(a)}_{m}!$ is equal to $R(\nu,N)$.
\end{thm}
The proof is due to the inclusion-exclusion 
principle and an explicit evaluation of the 
M\"obius inversion formula similar to \cite{KN1,KN2}.

\section{$R(\nu,N)$ and $Q$-system}\label{sec:3}
So much for the Bethe equation, 
we now turn to the $Q$-system.
For $a,b \in I_\sigma$ and $m,k \in {\mathbb Z}$, set 
\begin{equation*}
G_{am,bk}=
\begin{cases}
-\frac{1}{\epsilon_b}
A_{ba}\delta_{m, k}
& r>1 \\
-A_{ba}
(\delta_{m,2k-1}+2\delta_{m,2k}+\delta_{m,2k+1})
&
d_b/d_a=2\\
-A_{ba}
(\delta_{m,3k-2}+2\delta_{m,3k-1}+3\delta_{m,3k}&
d_b/d_a=3\\
\qquad\qquad\qquad
+2\delta_{m,3k+1}+\delta_{m,3k+2})&\\
-A_{ab}\delta_{d_a m, d_b k}
& \text{otherwise}.
\end{cases}
\end{equation*}

Let $\alpha_a$ and $\Lambda_a$ ($a \in I_\sigma$)
be the simple roots and the fundamental weights
of ${\mathfrak g}_0$. We set
\begin{equation*}
x_a=e^{\epsilon_a\Lambda_a},\quad
\quad y_a=e^{-\alpha_a}, 
\end{equation*}
which are related as 
\begin{align}
\label{eq:yx4}
y_a=\prod_{b=1}^n x_b^{-A_{ba}/\epsilon_b}.
\end{align}

\begin{defn}\label{defn:q-system}
The system of equations ($Q_{0}^{(a)}(y) =1$)
\begin{equation}
\label{q-system}
(Q^{(a)}_{m}(y))^2=
Q_{m+1}^{(a)}(y) Q_{m-1}^{(a)}(y)
+
y_a^{m}
(Q^{(a)}_{m}(y))^2
\prod_{\scriptstyle
(b,k)\in H}
(Q^{(b)}_{m}(y))^{G_{am,bk}}
\end{equation}
for a family $(Q^{(a)}_{m}(y))_{(a,m)\in H}$  of 
power series of $y=(y_{a})_{a=1}^{n}$
with unit constant terms is called the $Q$-system.
\end{defn}
The factor $y_a^m$ in the RHS is absorbed away if 
(\ref{q-system}) is written in terms of 
the combination $x_a^mQ^{(a)}_m(y)$.
The resulting form of the $Q$-system  has originally appeared in \cite{KR}
($A^{(1)}_n, B^{(1)}_n, C^{(1)}_n, D^{(1)}_n$), 
\cite{K3} ($E^{(1)}_{6,7,8},F^{(1)}_4, G^{(1)}_2$) and 
\cite{HKOTT} (twisted case).
Definition \ref{defn:q-system} corresponds to an {\em infinite} 
$Q$-system in the terminology of \cite{KNT}.
Its solution is not unique in general.
Following \cite{KNT} we introduce 
\begin{defn}
A solution of (\ref{q-system}) is {\em canonical}
 if the limit $\lim_{m\to \infty} Q^{(a)}_{m}(y)$
exists in the ring  $\mathbb{C}[[y]]$ 
of formal power series of $y=(y_{a})_{a=1}^{n}$ 
 with the standard topology. 
\end{defn}
\begin{thm}(\cite{KNT})\label{th:moohantai}
There exists a unique canonical solution $(\bQ^{(a)}_m(y))_{(a,m) \in H}$ 
of the $Q$-system (\ref{q-system}). 
Moreover, for any $\nu \in \mathcal{N}$, it admits the 
formula:
\[
\prod_{(a,m) \in H}(\bQ^{(a)}_m(y))^{\nu^{(a)}_m} = R^\nu(y), 
\]
where the power series $R^\nu(y)$ is defined by 
\begin{equation*}
R^\nu (y)=
\sum_{N\in \mathcal{N}} R(\nu,N) \prod_{a=1}^n
y_a^{\sum_{m=0}^\infty m N^{(a)}_{m}}
\end{equation*}
in terms of the integer $R(\nu,N)$ in (\ref{eq:Rdef}).
\end{thm}
In the proof of the theorem \cite{KNT}, the expression $R(\nu,N)$ emerges 
{}from a general argument on the $Q$-system, which 
is independent of the Bethe equation.
Our main finding in this paper is that it coincides with 
the number of off-diagonal solutions to the SCE
obtained in Theorem \ref{thm:Rnds}.

Let us state the consequence of this fact 
in the light of the Kirillov-Reshetikhin conjecture.
Let $\mbox{ch}^{(a)}_m(x)$ denote the Laurent polynomial of $x=(x_a)_{a=1}^n$
representing the ${\mathfrak g}_0$-character of the KR module
$W^{(a)}_m(\zeta)$.
Then, ${\mathcal Q}^{(a)}_m(y):=x_a^{-m}\mbox{ch}^{(a)}_m(x)|_{x=x(y)}$,
where $x(y)$ is the inverse map of  (\ref{eq:yx4}),
is a polynomial of $y=(y_a)_{a=1}^n$
with the unit constant term.
We call ${\mathcal Q}^{(a)}_m(y)$ the 
{\em normalized ${\mathfrak g}_0$-character}
of $W^{(a)}_m(\zeta)$.
The normalized character of 
the ${\mathfrak g}_0$-module 
$W^\nu$ in (\ref{eq:qspace}) is given by
\begin{equation*}
{\mathcal Q}^\nu(y)=\prod_{(a,m)\in H} 
({\mathcal Q}^{(a)}_{m}(y))^{\nu^{(a)}_{m}}.
\end{equation*}
The Kirillov-Reshetikhin conjecture [KR] is formulated in [KNT] as
\begin{conj}\label{conj:qsys}
${\mathcal Q}^{(a)}_{m}(y) = \bQ^{(a)}_m(y)$ for any 
$(a,m)\in H$. 
%
\end{conj}
%
%
Combining 
Theorem \ref{th:moohantai} and  Conjecture \ref{conj:qsys}, 
we relate the weight multiplicity in the tensor product of KR modules 
to the number of off-diagonal solutions to the SCE:
\begin{cor}[Formal completeness of the Bethe ansatz at $q=0$]
\label{thm:main3}
Under Conjecture \ref{conj:qsys} one has 
\begin{equation*}
{\mathcal Q}^\nu(y) = R^\nu(y). 
\end{equation*}
\end{cor}
%
Conjecture \ref{conj:qsys} implies that 
$(\prod_{(a,m)\in H} x_a^{m\nu^{(a)}_m})R^\nu(y(x))$ 
is a Laurent polynomial invariant 
under the Weyl group of ${\mathfrak g}_0$.
In fact canonical solutions 
have also been obtained as linear combinations 
of characters of irreducible 
finite dimensional ${\mathfrak g}_0$-modules for 
$X^{(r)}_N=A^{(1)}_{n},
 B^{(1)}_{n},C^{(1)}_{n},D^{(1)}_{n}$ \cite{KR,HKOTY}, and 
 for $X^{(r)}_N =A^{(2)}_{2n-1},A^{(2)}_{2n},D^{(2)}_{n+1},D^{(3)}_{4}$ 
 \cite{HKOTT}. 
For the current status of Conjecture \ref{conj:qsys}, 
see section 5.7 of \cite{KNT}.


\begin{thebibliography}{i}

\bibitem [B]{B}
R.\ J. \ Baxter,
\textit{Exactly solved models in statistical mechanics}, Academic Press, 
London (1982).

\bibitem [Be]{Be}
H. \ A. \ Bethe,
\textit{Zur Theorie der Metalle, I. Eigenwerte und
Eigenfunktionen der linearen Atomkette},
Z. Physik {\bf 71} (1931) 205--231.

\bibitem [C]{C} V.\ Chari,
{\em On the fermionic formula and the Kirillov-Reshetikhin
 conjecture,}
math.QA/0006090.


\bibitem [CP1]{CP1} V.\ Chari and A.\ Pressley,
\textit{ Quantum affine algebras and their
representations},
Canadian Math.\ Soc.\ Conf.\ Proc.\ {\bf 16}
(1995) 59--78.

\bibitem [CP2]{CP2} V.\ Chari and A.\ Pressley,
{\em Twisted Quantum affine algebras,}
Commun.\ Math.\ Phys.\ {\bf 196} (1998) 461--476.

\bibitem [HKOTY]{HKOTY} G.\ Hatayama,  A.\ Kuniba,
M.\ Okado, T.\ Takagi and Y.\ Yamada, \textit{Remarks on fermionic
 formula},
Contemporary Math.\ {\bf 248} (1999) 243--291.

\bibitem[HKOTT]{HKOTT}
G. Hatayama, A. Kuniba, M. Okado , T. Takagi and Z. Tsuboi, 
\textit{Paths, Crystals and Fermionic Formula}, 
math.QA/0102113.

\bibitem [Kac]{Kac} V.~G.~Kac, \textit{Infinite dimensional Lie
algebras}, 3rd edition, Cambridge Univ. Press, Cambridge (1990).

\bibitem [K1]{K1} A.\ N.\ Kirillov,
\textit{Combinatorial identities and completeness
 of states for  the  Heisenberg magnet},
J.\ Sov.\ Math.\ {\bf 30} (1985) 2298--3310.


\bibitem [K2]{K2} A. N. Kirillov, \textit{Completeness of states of
the generalized Heisenberg magnet}, J.\ Sov.\ Math. \ {\bf 36} (1987)
115--128.


\bibitem [K3]{K3} A. N. Kirillov,
\textit{Identities for the Rogers dilogarithm function
connected with simple Lie algebras},
 J.\ Sov. \ Math. \ {\bf 47} (1989) 2450--2459.


\bibitem[KR]{KR}
A.\ N.\ Kirillov and N.\ Yu.\ Reshetikhin,
\textit{Representations of Yangians and multiplicity of
occurrence of the irreducible components of the
tensor product of representations of simple Lie algebras},
J.\ Sov.\ Math.\ {\bf 52} (1990) 3156--3164.

\bibitem[KN1]{KN1}
A.\ Kuniba and T.\ Nakanishi,
\textit{The Bethe equation at $q=0$,
the M\"obius inversion formula, and weight multiplicities:
I.\ The $\mathfrak{sl}(2)$ case},
Prog.\ in Math.\ 191 (2000) 185--216.

\bibitem[KN2]{KN2}
A.\ Kuniba and T.\ Nakanishi,
\textit{The Bethe equation at $q=0$,
the M\"obius inversion formula, and weight multiplicities:
II.\ The $X_n$ case}, 
math.QA/0008047, J. Alg. in press.

\bibitem[KNT]{KNT}
A.\ Kuniba, T.\ Nakanishi and Z.\ Tsuboi,
\textit{The canonical solutions of the $Q$-systems and the Kirillov-Reshetikhin 
conjecture}, math.QA/0105145.




\bibitem [OW]{OW}
E.\ Ogievetsky and P.\ Wiegmann,
\textit{Factorized $S$-matrix and the Bethe ansatz for simple Lie
groups},
Phys.\ Lett.\ B {\bf 168} (1986) 360--366.



\bibitem [RW]{RW}
N.\ Yu.\ Reshetikhin and P.\ Wiegmann,
\textit{Towards the classification of completely integrable
quantum field theories (the Bethe ansatz associated with
Dynkin diagrams and their automorphisms)},
Phys.\ Lett.\ B {\bf 189} (1987) 125--131.

\end{thebibliography}
\end{document}